\documentclass[12pt]{amsart}

\usepackage{latexsym}

\def\conv{{\rm Conv \,}}
\def\R{{\bf R \,}}

\def\inter{{\rm int \,}}

\def\beq{\begin{equation}}
\def\eeq{\end{equation}}

\def\cal{\bf}

\newtheorem{theorem}{Theorem}
\newtheorem{lemma}{Lemma}
\newtheorem{corollary}{Corollary}

\newtheorem{assumption}{Assumption}
\newtheorem{definition}{Definition}

\newtheorem{example}{Example}
\newtheorem{remark}{Remark}
\newcommand{\half}{\mbox{${1 \over 2}$}}
\newcommand {\ve}{{\rm{vec}}}

\def\ba{\begin{array}}
\def\ea{\end{array}}
\def\beann{\begin{eqnarray*}}
\def\eeann{\end{eqnarray*}}
\def\bea{\begin{eqnarray}}
\def\eea{\end{eqnarray}}

\def\BT{\begin{theorem}}
\def\ET{\end{theorem}}
\def\BL{\begin{lemma}}
\def\EL{\end{lemma}}
\def\BC{\begin{corollary}}
\def\EC{\end{corollary}}
\def\BE{\begin{example}}
\def\EE{\end{example}}
\def\BD{\begin{definition}}
\def\ED{\end{definition}}
\def\BR{\begin{remark}}
\def\ER{\end{remark}}
\def\BAS{\begin{assumption}}
\def\EAS{\end{assumption}}

\def\la{\langle}
\def\ra{\rangle}

\begin{document}

\thispagestyle{empty}

\date{\today}

\author{
Vincent D. Blondel}
\address{Division of Applied Mathematics,
Universit\'e catholique de Louvain, 4 avenue Georges Lemaitre,
B-1348 Louvain-la-Neuve, Belgium}
\email{blondel@inma.ucl.ac.be}

\author{
Yurii Nesterov}
\address{Center for Operations Research and Econometrics (CORE),
Catholic University of Louvain (UCL), 34 voie du Roman Pays, 1348
Louvain-la-Neuve, Belgium} \email{nesterov@core.ucl.ac.be}

\title{Computationally efficient approximations of the joint spectral radius}


\maketitle

Abstract. {The joint spectral radius of a set of matrices is a
measure of the maximal asymptotic growth rate that can be obtained
by forming long products of matrices taken from the set.  This
quantity appears in a number of application contexts but is
notoriously difficult to compute and to approximate. We introduce
in this paper a procedure for approximating the joint spectral
radius of a finite set of matrices with arbitrary high accuracy.
Our approximation procedure is polynomial in the size of the
matrices once the number of matrices and the desired accuracy are
fixed.

For the special case of matrices  with non-negative entries we
give elementary proofs of  simple inequalities that we then use to
obtain approximations of arbitrary high accuracy. From these
inequalities it follows that the spectral radius of matrices with
non-negative entries is given by the simple expression
$$\rho(A_1, \ldots,  A_m)= \lim_{k
\rightarrow \infty} \rho^{1/k} ( A_1^{\otimes k} + \cdots +
A_m^{\otimes k})$$ where it is somewhat surprising to notice that
the right hand side does not directly involve any mixed product
between the matrices ($A^{\otimes k}$ denotes the $k$-th Kronecker
power of $A$).

For matrices with arbitrary entries (not necessarily non-negative)
we introduce an approximation procedure  based on semi-definite
liftings that can be implemented in a recursive way. For two
matrices, even the first step of the procedure gives an
approximation whose relative accuracy is at least ${1 /
\sqrt{2}}$, that is, more than $70\%$. The subsequent steps
improve the accuracy but also increase the dimension of the
auxiliary problems from which the approximation can be found.

Our approximation procedures  provide approximations of relative
accuracy $1-\epsilon$ in time polynomial in $n^{(\ln m)
/\epsilon}$, where $m$ is the number of matrices and $n$ is their
size.  These bounds are close from optimality since we show that,
unless P=NP, no approximation algorithm is possible that provides
a relative accuracy of $1-\epsilon$ and runs in time polynomial in
$n$ and $1/\epsilon$.

As a by-product of our results we prove that a widely used
approximation of the joint spectral radius based on common
quadratic Lyapunov functions (or on ellipsoid norms)  has relative
accuracy $1/\sqrt m$, where $m$ is the number of
matrices.}\footnote{This paper presents research results of the
Belgian Program on Interuniversity Poles of Attraction initiated
by the Belgian State, Prime Minister's Office, Science Policy
Programming. The scientific responsibility is assumed by the
authors.}

\section{Introduction}
\setcounter{equation}{0}

Let $\| \cdot \|$ be a matrix norm. The {\it spectral radius} of
the real matrix $A$ is defined by \beq \label{e0} \rho(A)= \lim_{k
\rightarrow +\infty} \|A^k\|^{1/k}.\eeq The spectral radius of a
matrix does not depend on the chosen matrix norm and it is also
equal to $\rho(A)= \max \{ |\lambda|: \lambda \mbox{ is an
eigenvalue of } A\}$ (see for example \cite[Corollary
5.6.14]{HJ}). The definition of spectral radius can be extended to
\emph{sets} of matrices in a  natural way. The \emph{joint
spectral radius} of a  set of matrices is a quantity introduced by
Rota and Strang in the early 60's  that measures the maximal
asymptotic growth rate that can be obtained by forming long
products of matrices; see \cite{rota}. The general definition is
for arbitrary sets of matrices but we shall consider here only
finite sets. Let $\{A_1, \ldots, A_m\}$ be some set of real
matrices. To the finite sequence $\sigma=(\sigma_1,  \sigma_2,
\ldots, \sigma_k) \in \{1,\dots,m\}^k$ we associate the
corresponding matrix product
$$A_\sigma= A_{\sigma_k} \cdots A_{\sigma_2} A_{\sigma_1}.$$ With
this notation, the joint spectral radius is defined by
\beq\label{def-rho} \rho(A_1, \dots , A_m) =  \limsup \limits_{k
\to +\infty}\; \max\limits_{\sigma \in \{1,\dots,m\}^k} \left\|
A_\sigma \right\|^{1/k}. \eeq

As for the single-matrix case, the joint spectral radius does not
depend on the matrix norm used. To see this, remember that any two
matrix norms $\|.\|_1$ and $\|.\|_2$ are related by $\alpha
\|A\|_1 \leq \|A\|_2 \leq \beta \|A\|_1$ for some $0 < \alpha <
\beta$. For any product $\sigma \in \{1,\dots,m\}^k$ one has
$\alpha^{1/k} \|A_\sigma\|_1^{1/k} \leq \|A_\sigma\|_2^{1/k} \leq
\beta^{1/k} \|A_\sigma\|_1^{1/k}$ and by letting $k$ tend to
infinity we conclude that the joint spectral radius is well
defined independently of the matrix norm used.

A definition analogous to (\ref{def-rho}) is possible by replacing
the norm appearing in the definition by a spectral radius. The
quantity defined in this way is the \emph{generalized spectral
radius} introduced in \cite{DauSe9222}. In  \cite{berger-wang},
the joint and generalized spectral radii of finite (or bounded)
sets of matrices are proved to be equal (see also \cite[Theorem
1]{elsner} for an elementary proof); thus reinforcing the status
of the joint spectral radius as a legitimate generalization of the
spectral radius of a single matrix. In the sequel we shall deal
only with the spectral radius defined with a norm, as in
(\ref{def-rho}).


Since its introduction in the 60's the joint spectral radius has
appeared in a number of different contexts; see, e.g.,
\cite{strang} or \cite{Dau2} for recent short
surveys\footnote{Google returns 625 entries upon entry of the
query ``joint spectral radius".}. Let us illustrate one
application in a dynamical system context. Consider the simple
discrete-time linear inclusion
$$x_{k+1} \in \{A_1, \ldots, A_m\} \; x_k, \quad x_0 \in \R^n$$ in
which at each step  a particular linear transformation is chosen
from of a finite number of possible choices. The maximal
asymptotic rate of growth of the trajectories associated to such a
system is given by the joint spectral radius $\rho(A_1, \ldots,
A_m)$. (In the discrete linear inclusion literature, the logarithm
of the joint spectral radius is sometimes called {\it Lyapunov
indicator}, see for example \cite{barabanov}.) In particular, all
possible trajectories will converge to the origin if and only if
$\rho(A_1, \ldots, A_m)<1$. The joint spectral radius can thus be
associated with the stability properties of time-varying linear
systems in the worst case over all possible time variations. It
also occur in the context of ``asynchronous'' \cite{tsitsiklis2}
or ``desynchronised'' \cite{Koz} systems. Besides
systems-theoretic interpretations, the concept is pervasive in
many areas of applied mathematics such as in wavelets
\cite{DauSe9222}, iterated function systems, random walks,
fractals, numerical solutions to ordinary differential equations
\cite{zen}, discrete-event systems \cite{bl}, interpolation
\cite{yu}, and coding theory \cite{moison}. \\


Despite its natural interpretation, the joint spectral radius  is
difficult to compute. Questions related to its computability and
to the existence of efficient approximation algorithms have been
posed more than a decade ago (see \cite{Topen} and \cite{LW95}).
In principle, the spectral radius can be approximated to any
desired accuracy
 by computing converging
sequences of upper and lower bounds. The following bounds, proved
in \cite{LW95}, \beq \label{e5} \max\limits_{\sigma \in
\{1,\dots,m\}^k} \rho (A_\sigma)^{1/k} \leq \rho(A_1, \ldots, A_m)
\leq \max\limits_{\sigma \in \{1,\dots,m\}^k} \| A_\sigma\|^{1/k}
\eeq can be evaluated for increasing values of $k$ and lead to
arbitrary accurate approximations of $\rho$ (see, e.g.,
\cite{Dau2} or \cite{grip}). Such approximation algorithms can in
turn be used in procedures that decide, after finitely many steps,
whether $\rho
>1$ or $\rho <1$ (such procedures are given, e.g.,  by Brayton and
Tong \cite{brayton} in a system theory context and by Barabanov
\cite{barabanov} in the context of discrete linear inclusions).
These procedures may however not terminate when $\rho$ happens to
be equal to 1 and the existence of algorithms for computing
arbitrarily precise approximations of $\rho$ does therefore not
rule out the possibility that the decision problem ``${\rho}<1$?''
is undecidable. It is so far unknown whether this is the case or
not (see \cite{LW95} for a discussion of this issue and for a
description of its connection
 with the finiteness conjecture that has since then been proved to be false; see
 \cite{bousch} and \cite{BTV}). A negative result in this direction is given in
 \cite{BTb} where it is proved that the related problem ``$\rho \leq 1$?" is algorithmically undecidable.\\

Approximations of the joint spectral radius that are directly
based on the inequalities (\ref{e5}) are expensive to compute. In
\cite{Maesumi}, the exponential number of products that appear in
the direct and naive computation of the bounds in (\ref{e5}) is
reduced by avoiding duplicate computation of cyclic permutations;
the total number of product to consider remains however
exponential. In addition to this, there are to this date no known
theoretical guarantees for the rate of convergence of the bounds
appearing in (\ref{e5}). Approximations of arbitrary degree of
accuracy can be computed, but at a price that may happen to be
prohibitive.\\

There are in fact intrinsic limitations for the rate at which the
joint spectral radius can be approximated. Let us say that the
value $\xi$ approximates the value $\rho$ with {\em relative
accuracy} $\mu \in [0,1]$ (or $100 \; \mu \; \%$), if $\mu \; \xi
\leq \rho \leq \xi$. By using a small adaptation of a proof
appearing in  \cite{BT} we show that, unless P=NP, there is no
algorithm that can compute the joint spectral radius of two
matrices with relative accuracy $1-\epsilon$ in time polynomial in
the size   of the matrices and in $1/\epsilon$ (see later for more
precise definitions). Despite this negative result it is still
conceivable, as pointed in \cite{Bsur}, that for any {\em fixed\/}
desired accuracy, there exists a polynomial time algorithm that
computes the joint spectral radius of the matrices with that
accuracy. We prove in this paper that this is indeed the case. In
Theorem \ref{t67}, we prove that the joint spectral radius of $m$
matrices of size $n$ can be approximated with relative accuracy
$1/\sqrt m$ by computing the spectral radius of a single matrix
whose size is less than $n^2$. This procedure can be applied in a
recursive way and in general we show how a relative accuracy of
$(1/\sqrt m)^{1/k}$ can be obtained by computing the spectral
radius of a single matrix of size less than $n^{2k}$. As an
illustration, the spectral radius of two matrices of size $n$ can
be computed with an accuracy of $70\%$ by computing the spectral
radius of a single matrix of size $n^2$, and an accuracy of 95\%
can be obtained for the joint spectral radius of three matrices by
computing the spectral radius of a single matrix of size $n^{11}$.
More generally, for any  number of matrices in the set and desired
relative accuracy we construct in Section \ref{s5} a single matrix
whose spectral radius approximates the joint spectral radius with
the desired accuracy. The approximation procedure runs in
polynomial time once the desired accuracy and the number of
matrices in the set are fixed. More precisely, our approximation
procedures provide approximations of relative accuracy
$1-\epsilon$ in time polynomial in $n^{(\ln m)/\epsilon}$, where
$m$ is the number of matrices and $n$ is their size. Notice that
$n^{(\ln m)/\epsilon} =e^{(\ln n \ln m)/\epsilon}=m^{(\ln
n)/\epsilon}$ and so our approximation procedure also runs in time
polynomial in the number of matrices once the desired accuracy and
the size of the matrices are fixed. These bounds are close from
optimality since we prove that, unless P=NP, no approximation
polynomial time
algorithm in $n$ and $1/\epsilon$ is possible.\\


 We now briefly
describe how our results are obtained and how the paper is
organized. We first notice in Section \ref{s2} that for matrices
with non-negative entries the joint spectral radius satisfies \beq
\label{iu9} \ba{c}\frac{1}{m} \ea \rho(A_1+ \cdots + A_m) \leq
\rho(A_1, \ldots,  A_m) \leq \rho(A_1+ \cdots + A_m).\eeq (The
left-hand side inequality is valid for arbitrary matrices; there
is no need to assume that the matrices have non-negative entries.)
Matrices with non-negative entries are exactly those matrices $A$
that are such that $A \;  \R_+^n \subseteq \R_+^n$ and we prove in
Theorem \ref{t1} that the inequalities (\ref{iu9}) are satisfied
not only for sets of matrices with non-negative entries but also
for sets of matrices that leave a proper cone invariant, i.e.,
matrices $A_i$ that are such that  $A_i \; K \subseteq K$ for some
proper cone $K$ and all $i$. Thus for these matrices the
inequalities (\ref{iu9}) provide a relative accuracy $1/m$. This
accuracy can be improved by considering Kronecker powers  of
matrices. In Section \ref{s3} we give an elementary proof that the
joint spectral radius of the $k$-th Kronecker powers of the
matrices in a set is equal to the $k$-th power of the joint
spectral radius of the set (Theorem \ref{t3}). Combining this with
the relation (\ref{iu9}) we prove that the approximation
$$\rho^{1/k}(A_1^{\otimes k} + \cdots
+A_m^{\otimes k})$$ has relative accuracy $1/m^{1/k}$. Some of the
arguments used to derive this result are valid only for matrices
that have  a common proper invariant cone. In Section \ref{s4} we
show how similar arguments can be used for arbitrary matrices. We
introduce a semi-definite lifting procedure that transforms a
linear operator acting $\R^n$ into a linear operator acting on the
space of symmetric $n \times n$ matrices. It appears that under
this transformation the joint spectral radius is simply squared.
Moreover, the lifting has the interesting feature that the
operator defined in this way leaves the cone of semi-definite
matrices invariant. This observation  leads to a simple
approximation procedure (Theorem \ref{t5}) for arbitrary matrices.
Implementation issues and numerical examples are provided in
Section \ref{s6}. Our approximation procedure provides a relative
accuracy $1-\epsilon$ in time polynomial in $n^{(\ln m)
/\epsilon}$. We prove in Section \ref{s7} (Theorem \ref{t6}) that
the same accuracy cannot be obtained in time polynomial in $n$ and
$1/\epsilon$ unless P=NP. The semi-definite lifting introduced in
Section \ref{s4} is a fundamental and powerful tool for stability
analysis. As an illustration of this, we prove in Section \ref{s8}
how to apply our results to analyze the quality of the so-called
ellipsoid approximation of the  joint spectral radius. We prove
that the ellipsoid approximation (a notion that is formally
defined in \cite{BTV} but that is implicitly present in a number
of earlier contributions on hybrid and time-varying systems, see,
e.g., \cite{gur} and \cite{ando} as well as some of the references
in \cite{LiberzonHespanhaMorse}) is an approximation of guaranteed
accuracy $1/\sqrt m$; a proof of this corollary can also be
extracted from Section 3 of \cite{gur}. For matrices of high
dimension this result significantly improves the earlier bound of
$1/\sqrt n$ proved in \cite{BTV} and in \cite{ando}. Finally, in a
last section, we discuss some of our results.

\section{Approximation for matrices leaving a cone invariant}\label{sc-Def}
\label{s2} \setcounter{equation}{0}

In this section, we consider sets of matrices that leave a proper
cone invariant and we show with elementary arguments how joint
spectral radius approximations of guaranteed accuracy can easily
be computed for this case. We start with a proof that the spectral
radius of a convex combination of matrices is always less or equal
to the joint spectral radius of the matrices. This result is valid
for all sets of matrices (there is no need to assume that the
matrices have non-negative entries) and is proved in \cite{BNT}
using the main result of \cite{Koz}. Here we present a direct and
elementary justification.

\BL For any set of matrices $\{ A_i: i=1, \ldots, m \}$ and any
$\alpha_i \geq 0$ satisfying $\sum_{i=1}^n \alpha_i=1$ we have:
\beq\label{eq-Bounds} \rho ( \sum\limits_{i=1}^m \alpha_i A_i)
\leq \rho(A_1, \dots , A_m). \eeq \EL

\proof Let us fix some $\alpha_i \geq 0$ with $\sum_i \alpha_i=1$
and an integer $k \geq 1$. Then
$$\|(\sum\limits_{i=1}^m
\alpha_{i} A_i )^k \| = \|\sum\limits_{\sigma \in \{1,\dots,m\}^k}
\alpha_\sigma A_\sigma\| \leq \sum\limits_{\sigma \in
\{1,\dots,m\}^k} \alpha_{\sigma} \|A_\sigma \| \leq
\max\limits_{\sigma \in \{1,\dots,m\}^k} \|A_\sigma\|$$ but then
also
$$\lim_{k\rightarrow \infty}\|(\sum\limits_{i=1}^m
\alpha_{i} A_i )^k \|^{1/k} = \limsup_{k\rightarrow
\infty}\|(\sum\limits_{i=1}^m \alpha_{i} A_i )^k \|^{1/k} \leq
\limsup_{k\rightarrow \infty} \max\limits_{\sigma \in
\{1,\dots,m\}^k} \|A_\sigma\|^{1/k}$$ and the result then follows
from the definitions of the spectral radius (\ref{e0}) and  of the
joint spectral radius (\ref{def-rho}).

 \qed

An immediate corollary is given by:

 \BC\label{cor-LowBound} $$\ba{c} {1
\over m} \ea \rho(\sum\limits_{i=1}^m A_i ) \leq \rho(A_1, \dots ,
A_m).$$ \qed \EC

The example $A_1=A, A_2=-A$ clearly shows that we cannot hope in
general to have $\rho(A_1, \dots , A_m) \leq \rho
(\sum\limits_{i=1}^m A_i)$. This inequality is nevertheless
satisfied when the matrices $A_i$ leave a proper cone invariant. A
{\em cone} in $\R^n$ is a subset $K \subseteq \R^n$ such that
$\lambda v \in K$ for all $\lambda \geq 0$ and $v \in K$. We say
that a cone $K$ is {\em proper} if it is closed, convex, has
nonempty interior, and contains no straight line. For example, the
set of vectors with non-negative entries $\R_+^n$ is a proper cone
but $\R^n$ itself is not (see \cite{rock} for more background on
cones and proper cones). We shall say that the matrices $A_i$ {\em
leave a proper cone invariant} if there exists a proper cone $K
\subseteq \R^n$ such that  $A_i \, K \subseteq K$ for all $i$. For
example, matrices with non-negative entries leave the proper cone
$\R_+^n$ invariant.

\BL\label{l5}
 Associated to any proper cone $K$ there is a matrix norm $\|
 \cdot \|_K$ that satisfies $\|A\|_K \leq
\|A+B\|_K$ for all matrices $A$ and $B$ that leave the cone $K$
invariant. \EL

\proof  Most usual matrix norms satisfy this property when
$K=\R_+^n$; we provide a construction for arbitrary proper cones
$K$.

Let $\| \cdot \|$ be some arbitrary vector norm. The dual vector
norm $\|\cdot \|_*$ is defined by $\|w\|_*=\max_{\|v\|=1}w^Tv.$
Assume that $K$ is a proper cone such that $A_i \; K \subseteq K$.
Define the dual of $K$ by
$$K^*=\{w \in \R^n: w^Tv \geq 0 \mbox{ for all } v \in K\}.$$ The
dual of a proper cone is again a proper cone, see \cite{rock}. Let
us now consider the quantity
$$\|A\|_K=\max\limits_{{v \in K, w \in K^*}  \atop {\|v\|=\|w\|_*=1}}
{w^TAv}.$$ It is easy to verify that because $K$ and its dual are
proper cones, the quantity $\|\cdot \|_K$  is indeed a matrix
norm. Moreover, by the definition of the norm,  matrices $A$ and
$B$ that satisfy $A K \subseteq K$ and $B K \subseteq K$, are such
that $\|A\|_K \leq \|A+B\|_K$.
 \qed

We may now state the result of this section.

 \BT\label{th-Positive}
 \label{t1}
 Let $\{ A_i: i=1, \ldots, m \}$ be a set
 of matrices that leave a proper cone invariant. Then \beq\label{eq-Positive} \ba{c} {1 \over m}
\rho(\sum\limits_{i=1}^m A_i) \leq \rho(A_1, \dots , A_m) \leq
\rho(\sum\limits_{i=1}^m A_i). \ea \eeq \ET

\proof The lower bound in (\ref{eq-Positive}) is already
established in Corollary \ref{cor-LowBound}. Assume that the
matrices in  $\{ A_i: i=1, \ldots, m \}$ leave the proper cone $K$
invariant and define the matrix norm $\|
 \cdot \|_K$ as in Lemma
\ref{l5}. Since $A_i K \subseteq K$, we also have $A_\sigma K
\subseteq K$ for all $\sigma \in \{1, \ldots, m\}^k$ but then also
$$\max_{\sigma \in \{1, \ldots, m\}^k}\|A_\sigma\|_K \leq \|\sum_{\sigma \in \{1, \ldots, m\}^k} A_\sigma\|_K = \|(\sum_{i=1}^m A_i)^k \|_K$$
and therefore
$$\rho(A_1, \ldots, A_m)=\limsup_{k\to\infty} \max_{\sigma \in \{1, \ldots, m\}^k}\|A_\sigma\|_K^{1/k}
\leq \limsup_{k\to\infty} \|(\sum_{i=1}^m A_i)^k
\|_K^{1/k}=\rho(\sum_{i=1}^m A_i).$$
 \qed

It is interesting to notice that the conclusion of Theorem
\ref{t1} does not directly involve the cone $K$. The existence of
a proper cone satisfying the hypotheses suffices to conclude and
the exact nature of the cone is irrelevant.

A matrix with non-negative entries leaves the cone $\R_+^n$
invariant and so we have the following corollary.

\BC\label{ex-Positive} Let the matrices $A_i$, $i=1,\dots,m$, have
non-negative entries. Then $$\ba{c}\frac{1}{m} \ea
\rho(\sum\limits_{i=1}^m A_i) \leq \rho(A_1, \dots , A_m) \leq
\rho (\sum\limits_{i=1}^m A_i).$$ \qed \EC


\section{Kronecker lifting}\label{sc-Def}
\label{s3}

We now describe a way to exploit Theorem \ref{t1} for obtaining
approximations of arbitrary high accuracy. Our approximations
involve Kronecker powers of matrices. The {\em Kronecker product}
of two matrices $A \in R^{p_1 \times q_1}$ and $B \in R^{p_2
\times q_2}$ is the ${p_1p_2 \times q_1q_2}$ matrix defined by
$$
A \otimes B = \left( \ba{ccc} A_{1,1} B & \dots &
A_{1,q_1} B \\
& \dots & \\
A_{p_1,1} B & \dots & A_{p_1,q_1} B \ea \right).
$$

We also define the {\em Kroneker power} of a matrix:
$$
A^{\otimes k} =\underbrace{A \otimes A \ldots A \otimes A}_{k
\mbox{ times}}.
$$
There is no need for parenthesis in this expression since the
Kronecker product of two matrices is an associative operation. Let
us prove some elementary properties.

\BL\label{lm-KronProp} \label{l2}  For matrices of appropriate
sizes we have:
\begin{enumerate}
\item $(A \otimes B)^T=A^T \otimes B^T$; \item $(A_1 \otimes B_1)
(A_2 \otimes B_2)=(A_1 A_2) \otimes (B_1 B_2)$; \item $(A^{\otimes
k})^T=(A^T)^{\otimes k}$; \item $A^{\otimes k} \, B^{\otimes k} =
(A B)^{\otimes k}$; \item Let $\|\cdot \|$ denote the spectral
matrix norm induced by the standard Euclidean vector norm. Then
$\|A^{\otimes k}\|=\|A\|^k$.
\end{enumerate} \EL

\proof

The first statement directly follows from the definition of the
Kronecker product. The second statement is quite standard; see,
e.g., Section 4.2 in \cite{HJ}. The third and fourth statement
follow, respectively, from repeated application of the first and
second statement. Finally, for proving the last statement, note
that by  using (3) and (4) we obtain

$$
\|A^{\otimes k}\|^2=\max\limits_{\| x \| = 1} x^T (A^{\otimes
k})^T A^{\otimes k} x= \max\limits_{\| x \| = 1} x^T
(A^T)^{\otimes k} A^{\otimes k} x= \max\limits_{\| x \| = 1} x^T
(A^T A)^{\otimes k} x.
$$

The matrix $(A^T A)^{\otimes k}$ is symmetric and  so the last
expression is also equal to the largest magnitude of the
eigenvalues of $(A^T A)^{\otimes k}$. For a matrix $B$ the
eigenvalues of the matrix $B^{\otimes k}$ are given by all
possible products of $k$ eigenvalues of $B$ (see Theorem 4.2.12 in
\cite{HJ}) and we therefore conclude $\max\limits_{\| x \| = 1}
x^T (A^T A)^{\otimes k} x = \|A\|^k$. \qed

We can now prove a useful identity for the spectral radius of
Kronecker powers of matrices.

\BT\label{th-Kron}\label{t2} Let $\{ A_i: i=1, \ldots, m \}$ be a
set of matrices and  $l \geq 1$. Then
$$\rho(A_1^{\otimes l}, \ldots, A_m^{\otimes l}) =
\rho^l(A_1,\ldots, A_m).$$ \ET

\proof Let $\sigma \in \{1, \ldots, m\}^l$. By using the
equalities (3) and (5) of Lemma \ref{l2} we obtain
$$\|A^{\otimes
l}_{\sigma_l}  \cdots  A^{\otimes l}_{\sigma_1}\| = \|
A_\sigma^{\otimes l} \| = \| A_\sigma \|^l.
$$
The result then follows from the definition of the joint spectral
radius (\ref{def-rho}). \qed

We now would like to combine Theorem \ref{th-Kron} with Theorem
\ref{t1}  in order to obtain approximations of arbitrary high
accuracy for matrices leaving a cone invariant. In order to do
this we need to consider Kronecker products of cones. For two sets
 $Q_1 \subseteq R^n$ and $Q_2 \subseteq R^m$ let us denote
$$
Q_1 \otimes Q_2 = \{ z = u \otimes v:\; u \in Q_1, \; v \in Q_2
\}.
$$
(In this definition, the vectors $u$ and $v$ are treated as column
matrices.) If $K_1$ and $K_2$ are proper cones then $K_1 \otimes
K_2$ does not need to be a proper cone. Indeed, consider for
example the proper cones $K_1=K_2=\R_+^2$ for which it is easy to
see that $K_1 \otimes K_2$ has empty interior in $\R^4$. We do
however have the following result.

\BL \label{l3} Let $K_1$ and $K_2$ be proper cones, then the cone
$K=\overline{\conv (K_1 \otimes K_2)}$ is a proper cone also. \EL

 \proof The cone $K$ is closed and convex by
definition. Assume that $\inter K = \emptyset$. Then there exists
$A \in R^{nm}\setminus \{ 0 \}$ such that
$$
\la A, z \ra = 0 \quad \forall z \in K \supset K_1 \otimes K_2.
$$
Note that the function $f(u,v) = \la A, u \otimes v \ra$ is a
bilinear form defined on $K_1 \times K_2$. Since $\inter (K_1
\times K_2) \neq \emptyset$, $f$ can vanish on this set only if $A
= 0$. That is a contradiction.

Let us prove now that $K$ contains no straight line. Assume that
this is not the case. Then $\inter K^* = \emptyset$, which implies
existence of $B \in R^{nm}\setminus \{0\}$ such that
\beq\label{eq-vanish} \la B, w \ra = 0 \quad \forall w \in K^*.
\eeq Since $\la u \otimes v, x \otimes y \ra = \la u, x \ra \cdot
\la v, y \ra$, we have $K_1^* \otimes K_2^* \subseteq K^*$. Since
$K_1^*$ and $K_2^*$ have nonempty interior, relation
(\ref{eq-vanish}) is impossible by the first part of the proof.
\qed

We are  ready to state the main result of this section.

\begin{theorem}
\label{t3}  Let $\{ A_i: i=1, \ldots, m \}$ be a set of matrices
leaving a proper cone invariant. Then
$$\ba{c} \frac{1}{m^{1/k}} \ea \rho^{1/k}(A_1^{\otimes k} + \cdots + A_m^{\otimes k}) \leq
\rho(A_1, \ldots, A_m)\leq \rho^{1/k}(A_1^{\otimes k} + \cdots
+A_m^{\otimes k}).$$ In particular, the joint spectral radius is
given by \beq \label{e31}\rho(A_1, \ldots, A_m)= \lim_{k
\rightarrow \infty} \rho^{1/k}(A_1^{\otimes k} + \cdots +
A_m^{\otimes k}).\eeq
\end{theorem}

\proof Assume that the matrices leave the proper cone $K$
invariant. If $A K \subseteq K$, then $(AK)^{\otimes k} \subseteq
K^{\otimes k}$ and by Lemma \ref{l2}, $A^{\otimes k} K^{\otimes k}
\subseteq K^{\otimes k}$. But then also $A^{\otimes k} \;
\overline{\conv K^{\otimes k}} \subseteq \; \overline{\conv
K^{\otimes k}}.$ By Lemma \ref{l3} the cone $\overline{\conv
K^{\otimes k}}$ is a proper cone. Thus the matrices $A_i^{\otimes
k}$ leave a proper cone  invariant and by Theorem \ref{t1} we have
$$\ba{c} \frac{1}{m} \ea \rho(A_1^{\otimes k} + \cdots + A_m^{\otimes k}) \leq
\rho(A_1^{\otimes k}, \ldots, A_m^{\otimes k})\leq
\rho(A_1^{\otimes k} + \cdots +A_m^{\otimes k}).$$ In order to
conclude it suffices to use the fact, proved in Theorem \ref{t2},
that
$$\rho(A_1^{\otimes k}, \ldots, A_m^{\otimes k}) =
\rho^k(A_1,\ldots, A_m).$$ \qed

It is somewhat surprising to notice that the right hand side in
(\ref{e31}) does not directly involve any mixed product between
the matrices.

An approximation of relative accuracy $(1/m)^{1/k}$ can thus be
obtained by computing the spectral radius of a single matrix of
dimension $n^k$. For pairs of matrices some of the relative
accuracies and the corresponding matrix sizes are as follows:

$$\ba{llllllllllll}
\mbox{Relative accuracy} & 0.707 & 0.840 & 0.917 & 0.957\\
\mbox{Matrix size} & n^{2} & n^4 & n^8 & n^{16} \ea$$


\section{Semi-definite lifting}\label{sc-Lift}
\label{s4} \label{s5} \setcounter{equation}{0}

The result in Theorem \ref{t3} provides an easy way to evaluate
with arbitrary accuracy the joint spectral radius of matrices that
leave a proper cone invariant. Unfortunately, matrices do not
always leave a proper cone invariant. We show in this section how
an invariant cone can always be created by semi-definite lifting.

Let $A \in \R^{n \times n}$ and consider the following linear
operator: \beq \label{e34} X \rightarrow A X A^T: \R^{n \times n}
\rightarrow \R^{n \times n}. \eeq

 A matrix representation for this linear operator can be
obtained by using the matrix-to-vector operator that develops a
matrix into a vector by taking its columns one by one. This
operator, denoted $\ve$, satisfies the elementary property  $\ve
(C X D)=(D^T \otimes C) \; \ve (X)$ (see Lemma 4.3.1 in
\cite{HJ2}). We therefore have $\ve (A X A^T)=(A \otimes A) \; \ve
X$ and a matrix representation of the linear operator (\ref{e34})
is thus simply given by $A \otimes A$.

Consider now the space of {\em symmetric} matrices ${\cal S}$.
This space is a subspace of $\R^{n \times n}$ of dimension
$n(n+1)/2$ and the operator \beq \label{e11} X \rightarrow A X
A^T: {\cal S} \rightarrow {\cal S}. \eeq is again a linear
operator. Matrix representations of this operator are of course
different from those  of the operator (\ref{e34}) since, in
particular, the spaces on which these operators operate have
different dimensions (we describe in the next section how to
construct a matrix representations for the operator on symmetric
matrices). In the next theorem we prove that even though the
operators $X \rightarrow A X A^T$ on ${\cal S}$ and on $\R^{n
\times n}$ are different, their joint spectral radius  are equal.

\BT\label{t4}\label{th-Square} Let $\{ A_i \in \R^{n \times n}:
i=1, \ldots, m \}$ and denote by  $M_A$ a matrix representation of
the linear operator $X \rightarrow A X A^T: {\cal S} \rightarrow
{\cal S}$. Then we have $$\label{eq-Square} \rho(M_{A_1}, \dots ,
M_{A_m}) = \rho(A_1 \otimes A_1, \dots , A_m \otimes
A_m)=\rho^2(A_1, \dots , A_m).$$ \ET

\proof The second equality is already proved in Theorem \ref{t2};
we only need to prove the first equality. Let $\| \cdot \|_F$
denote the Frobenius matrix norm (i.e., the sum of squares of all
entries) and consider the resulting induced operator norms for $X
\rightarrow A X A^T$ on $\R^{n \times n}$ and ${\cal S}$:

$$\sup_{X \in \R^{n \times n}, \|X\|_F=1} \|AXA^T\|_F \quad \mbox{
and } \quad \sup_{X \in {\cal S}, \|X\|_F=1} \|AXA^T\|_F.$$

We claim that these two operator norms are equal and that the
supremum is achieved for some symmetric matrix of rank one.
Indeed, using the
 matrix-to-vector operator $\ve$ and denoting by $\| \cdot \|$ the
usual vector Euclidean norm we get
\begin{eqnarray*}
\sup_{X \in \R^{n \times n}, \, \|X\|_F=1} \|AXA^T\|^2_F &=&
\sup_{x \in
\R^{n^2}, \|x\|=1} \|(A \otimes A) x\|^2\\
& = & \sup_{x \in \R^{n^2}, x^Tx =1} x^T (A \otimes A)^T (A
\otimes A) x\\
& = & \sup_{x \in \R^{n^2}, x^Tx =1} x^T (A^TA \otimes A^TA) x
\end{eqnarray*}

For deriving the last equality, we have used (2) in Lemma
\ref{l2}. The supremum in the last expression is achieved by any
eigenvector $x$ associated to an eigenvalue of largest magnitude
of $A^TA \otimes A^TA$. Let $(\lambda_i, v_i)$ $i=1, \ldots, n$ be
the set of eigenvalues/eigenvectors pairs associated to the
symmetric matrix $A^T A$ and let $\lambda_1$ be such that
$|\lambda_1|\geq \lambda_i$ for all $i$. The
eigenvalues/eigenvectors pairs of the matrix $A^TA \otimes A^TA$
are then given by $(\lambda_i \lambda_j, v_i \otimes v_j)$ for
$i,j=1, \ldots, n$. The vector $v_1 \otimes v_1$ is thus an
eigenvector of $A^TA \otimes A^TA$ of largest eigenvalue
magnitude. Since $\ve (v_k \, v_k^T)= v_k \otimes v_k$ we see that
the supremum of $\|AXA^T\|_F$ is achieved for the symmetric rank
one matrix $v_1 \, v_1^T$. From this we conclude that
$$\sup_{X \in \R^{n \times n}, \|X\|_F=1} \|AXA^T\|_F = \sup_{X \in {\cal S}, \|X\|_F=1} \|AXA^T\|_F$$
and thus also
$$\sup_{ \|x\|=1} \|(A \otimes A) x \| = \sup_{ \|x\|=1} \|(M_A) x \|.$$
Notice that $M_{A_\sigma}=(M_{A_i})_\sigma$ and so
$$\|(M_{A_i})_\sigma\| = \|M_{A_\sigma}\| = \| A_\sigma \otimes A_\sigma \| =
\|  (A_i \otimes A_i)_\sigma \|$$ and it then suffices to apply
the definition of the joint spectral radius to conclude the proof.

\qed

A matrix $A$ is {\em positive semi-definite} (which we denote by
$A \succeq 0$) if $ v^T A v \geq 0$ for all $v \in \R^n$. The set
of symmetric positive semi-definite matrices, denoted ${\cal
S}_+$, is a cone. It is not a proper cone of $\R^{n \times n}$
because it has empty interior in $\R^{n \times n}$; but it is a
proper cone of the set of symmetric matrices ${\cal S}$. If $X \in
{\cal S}_+$ then we clearly have $A X A^T \in {\cal S}_+$ and so
the operator $X \rightarrow A XA^T$ leaves the proper cone ${\cal
S}_+$ invariant. Combining this observation with Theorem \ref{t3}
and Theorem \ref{t4} we deduce:

\BT \label{t5}\label{t67} Let $\{ A_i \in \R^{n \times n}: i=1,
\ldots, m \}$ and denote by  $M_A$ a matrix representation of the
linear operator $X \rightarrow A X A^T: {\cal S} \rightarrow {\cal
S}$. Then \beq\label{eq-Bound1} \ba{c} {1 \over \sqrt{m}} \;
\rho^{1/2}(M_{A_1} + \ldots + M_{A_m}) \leq \rho(A_1, \dots , A_m)
\leq \rho^{1/2}(M_{A_1} + \ldots + M_{A_m}). \ea \eeq \ET

Note that for $m = 2$ the bounds (\ref{eq-Bound1}) are not so bad.
Indeed, ${1 /\sqrt{2}} > 0.7$ and so the relative accuracy of the
approximation is at least $70\%$.  The main interest of this
result resides however in the possibility of applying it
recursively. Starting from an initial matrix  $A \in \R^{n \times
n}$ we define a sequence of operators acting on spaces of
increasing dimensions. For $l = 1$ we define
$$
A^{[1]}(x) \equiv A x \quad x \in E_1 \equiv \R^n
$$
and for $l \geq 2$ we define the operators recursively
$$
\ba{rcl} A^{[l]}(x) & \equiv & M_{A^{[l-1]}}(x)  \quad x \in E_{l}
\equiv {\cal S}(E_{l-1}). \ea
$$
Assume for the simplicity of the presentation that $m=2$ and let
$A_1, A_2 \in \R^{n \times n}$. From Theorem \ref{th-Square} we
know that
$$
\rho\left(A_1^{[l]},A_2^{[l]} \right) = \rho^{2^l}(A_1,A_2)
$$
On the other hand, for any $l \geq 1$ the operators $A_1^{[l]}$
and $A_2^{[l]}$ leave the cone ${\cal S}_+(E_{l-1})$ invariant.
Therefore, in view of Theorem \ref{th-Positive},  we have
$$
\half \rho\left( A_1^{[l]}+A_2^{[l]}\right) \leq
\rho\left(A_1^{[l]},A_2^{[l]}\right) \leq
\rho\left(A_1^{[l]}+A_2^{[l]}\right).
$$
Combining these last two expressions, we get the following bounds:
$$
\ba{c} \left( \half \right)^{1/2^l}
\left[\rho\left(A_1^{[l]}+A_2^{[l]}\right)\right]^{1/2^l} \leq
\rho(A_1,A_2) \leq
\left[\rho\left(A_1^{[l]}+A_2^{[l]}\right)\right]^{1/2^l}.\ea
$$
Note that $$\ba{c}  1 - {\ln 2 \over 2^l} \leq \left( \half
\right)^{1/2^l} \leq 1 \ea $$ and thus the improvement in the
quality of our approximation is quite fast. Unfortunately, the
dimensions of the spaces $E_l$ are also growing fast; we have
$$
n_{l+1} = \half n_l(n_l+1)
$$
(where $n_i=\dim E_i$) and therefore asymptotically we have $n_l =
O\left(({n / 2})^{2^l}\right)$. Let us display, for a pair of
matrices, the relative accuracy of our approximation as a function
of the resulting dimension.
$$
\ba{|r|r|r|r|r|} \hline \mbox{Steps} & \mbox{Accuracy} & n = 2 & n
= 10 & n = 100
\\
\hline
1 & 0.707 & 3 & 55 & 5050\\
2 & 0.840 & 6 & 1540 & *\\
3 & 0.917 & 21 & 118570 & * \\
4 & 0.957 & 231 & * & * \\
5 & 0.978 & 26796 & * & * \\
\hline \ea
$$
We use the symbol $*$ to mark the cases for  which the dimension
of the auxiliary problem goes beyond the abilities of modern
computers.

\section{Numerical implementation and examples}
\label{s6} \setcounter{equation}{0}

The recursive  semi-definite lifting for obtaining approximations
of the joint spectral radius of increasing accuracy may be
difficult to implement because the matrix $M_A$ is not easy to
express in terms of the matrix $A$. Consider for example the case
of $2 \times 2$ matrices. Let
$$A= \left( \ba{cc} a_{11} & a_{12}\\ a_{21} & a_{22} \ea
\right),$$ then one easily compute  $$M_A=\left( \ba{ccc} a_{11}^2
& 2 a_{11}
a_{12} &  a_{12}^2\\
a_{11}a_{21} & a_{11} a_{22} + a_{12} a_{21} & a_{12} a_{22}
\\
a_{21}^2 & 2 a_{21} a_{22} & a_{22}^2 \ea \right).$$ In general,
the matrix $M_A$ can be expressed as follows: \beq \label{e13}
M_A=(Q^T Q)^{-1} Q^T P^{-1} (A \otimes A) P Q\eeq where $P$ is a
particular permutation matrix of size $n^2$ and $Q$ is given by
$$Q=\left( \ba{cc} I_n & 0 \\ 0 & I_{\frac{n(n-1)}{2}}\\ 0 &
I_{\frac{n(n-1)}{2}} \ea \right)$$ (the matrix $I_k$ is the
identity matrix of size $k$). One difficulty with the expression
(\ref{e13})  is that the permutation matrix $P$ is tedious to
construct; it essentially relates the column decomposition of a
matrix with its diagonal decomposition. The permutation matrix can
of course be constructed but the implementation of this
construction is somewhat cumbersome. A way to bypass this
difficulty can be achieved by using the semi-definite lifting only
once, and then apply only Kronecker liftings. The semi-definite
lifting can be performed either at the beginning or at the end of
the process, leading to the bounds \beq \label{e43}
\ba{c}\frac{1}{m}\ea \rho(M^{\otimes l}_{A_1}+\cdots + M^{\otimes
l}_{A_m}) \leq \rho^{2l}(A_1, \ldots, A_m) \leq \rho(M^{\otimes
l}_{A_1}+\cdots + M^{\otimes l}_{A_m}) \eeq and \beq \label{e44}
\ba{c}\frac{1}{m}\ea \rho(M_{A_1^{\otimes l}}+\cdots +
M_{A_m^{\otimes l}}) \leq \rho^{2l}(A_1, \ldots, A_m) \leq
\rho(M_{A_1^{\otimes l}}+\cdots + M_{A_m^{\otimes l}}). \eeq
 The validity of these inequalities results from the combination of Theorem \ref{t1},
 Theorem \ref{t2}, Lemma \ref{l3} and Theorem \ref{t4}.
 The expressions (\ref{e43}) and (\ref{e44}) both provide a relative accuracy of $(1/m)^{1/(2l)}$
  but involve matrices of
different size. The matrices $M^{\otimes l}_{A_i}$ have size
$(n(n+1)/2)^l$ whereas the matrices $M_{A_i^{\otimes l}}$ have
size $n^l(n^l+1)/2 > (n(n+1)/2)^l$. In addition to this, the
matrices $M^{\otimes l}_{A_i}$ are  easier to compute than
$M_{A_i^{\otimes l}}$ since they necessitate the evaluation of
semi-definite liftings of smaller matrices.

Of course, other combinations of Kronecker and semi-definite
liftings are also possible. Whenever a $k$-th Kronecker power is
used, the current relative accuracy is $k$-th rooted, and whenever
a semi-definite lifting is used, the relative accuracy is square
rooted. Moreover, at least one semi-definite lifting is needed if
the original matrices do not leave a proper cone invariant.

For the numerical implementation of these approximations notice
also that the spectral radius of $M_{A_1}+ \cdots + M_{A_m}$ is
the spectral radius of the linear operator
$$B: X \rightarrow A_1XA_1^T + \cdots + A_m X A_m^T: {\cal S} \rightarrow {\cal S}$$
and it can therefore be found by standard linear algebra
techniques, or by solving the following linear inequality
optimization problem \beq\label{prob-SDP} \ba{c} \inf\limits_{X,
\tau}\;\left\{ \tau:\; \tau X \succeq A_1 X A_1^T + \cdots + A_m X
A_2^T,\; X \succ 0 \right\} \ea \eeq for which efficient convex
optimization algorithms are available (see, e.g., \cite{boyd})
that are implemented in standard commercial softwares such as the
Matlab ``LMI Control Toolbox".

\section{Computational complexity analysis}
\label{s7} \setcounter{equation}{0}

The joint spectral radius can be approximated to arbitrary
accuracy. It is proved in \cite{BT} that, unless P=NP,
approximating algorithms of relative accuracy $1-\epsilon$ for
pairs of matrices cannot run in time polynomial in the size of the
matrices and in $\ln(1/\epsilon)$. It was later noticed (see the
note 9 on page 1260 of \cite{Bsur}) that a careful examination of
the proof in \cite{BT} shows that, unless P=NP,  there are no
approximation algorithm of relative accuracy $1-\epsilon$ that run
in time polynomial in the size of the matrices and in
$1/\epsilon$. This is however not proved in \cite{Bsur}. In this
section we prove this result by extracting the essential part of
the proof in \cite{BT}.

We proceed by reduction from the classical NP-complete
satisfiability problem 3-SAT. This problem is defined as follows.
Consider a set $\{x_1,\ldots x_n\}$ of Boolean variables. A
\emph{literal} is either a Boolean variable $x_i$, or its negation
${\bar x_i}$. A \emph{three-literal clause} is a disjunction of
three literals (e.g., $x_3 \vee \bar x_5 \vee \bar x_6$). The
3-SAT problem is the problem of determining, for a given
collection of three-literal clauses, if there exists a truth
assignment for the variables that simultaneously satisfies all
clauses. This problem is known to be NP-complete; finding a
polynomial time algorithm for the problem is exactly equivalent to
proving that P=NP.

The complexity result proved in \cite{BT} is based on the
construction of a pair of  matrices which we briefly outline (the
complete construction can easily be extracted from \cite{BT}).
From a given instance of 3-SAT with $n$ variables and $p$ clauses
two square matrices $A_1, A_2$ are constructed in  polynomial
time. The matrices have their entries in $\{0, 1\}$, are of size
$(n+1)(p+1)$, and have a joint spectral radius that satisfies:
$$\ba{rcll} \rho(A_1, A_2) & \geq & p^{1/(n+2)} & \mbox{ if the
instance is satisfiable}\\
  \rho(A_1, A_2) & \leq & (p-1)^{1/(n+2)} & \mbox{ if the instance is not satisfiable.} \ea
 $$
Any algorithm of relative accuracy $(1-1/p)^{1/(n+2)}$ allows to
make the  distinction between these two cases. Moreover, since
$$\left(1-\frac{1}{p}\right)^{1/(n+2)} \leq 1-\frac{1}{p(n+2)} \leq 1$$
it is clear that an approximation algorithm of relative accuracy
$1-{1}/{p(n+2)}$ also allows to decide 3-SAT. Since the size of
the matrices in the construction are given by $(n+1)(p+1)$, we
deduce:

\BT \label{t6} Unless P=NP, the problem of approximating the joint
spectral radius of two square matrices with $\{0, 1\}$ entries and
with relative accuracy $1-\epsilon$ cannot be obtained in time
polynomial in the size of the matrices and in $1/\epsilon$. \ET

From Theorem \ref{t3} we know that the spectral radius of two
matrices with nonnegative entries is given by $\rho(A_1, A_2)=
\lim_{k \rightarrow \infty} \rho^{1/k}(A_1^{\otimes k} +
A_2^{\otimes k})$ and so the quantity $\lim_{k \rightarrow \infty}
\rho^{1/k}(A_1^{\otimes k} + A_2^{\otimes k})$ is NP-hard to
approximate in the sense given in Theorem \ref{t6}. Theorem
\ref{t3} provides also a rate of convergence for the
approximations $\rho^{1/k}(A_1^{\otimes k} + A_2^{\otimes k})$.
From this we may prove the following complexity result.

\BT The problem of determining, for a given integer $k\geq 0$ and
for a given pair of matrices $A_1, A_2$ with nonnegative entries,
if $$\rho(A_1^{\otimes k}+ A_2^{\otimes k})<1$$ is a problem that
is NP-hard. \ET

\proof

We proceed by reduction from 3-SAT. From a given instance of 3-SAT
with $n$ variables and $p$ clauses we use the construction given
in \cite{BT} to construct two square matrices $B_1, B_2$ that have
their entries in $\{0, 1\}$, are of size $(n+1)(p+1)$, and have a
joint spectral radius that satisfies:
$$\ba{rcll} \rho(B_1, B_2) & \geq & p^{1/(n+2)} & \mbox{ if the
instance is satisfiable}\\
  \rho(B_1, B_2) & \leq & (p-1)^{1/(n+2)} & \mbox{ if the instance is not satisfiable.} \ea
 $$
 We then choose $\alpha
\in \R$ with $0 < \alpha <1$, $r \in {\bf Q}$ and $k\geq 0$ such
that $\alpha < (1/2)^{1/k}$ and \beq \label{e61}
(p-1)^{\frac{1}{n+2}} < \frac{1}{\alpha} (p-1)^{\frac{1}{n+2}} <
r<  \alpha q^{\frac{1}{n+2}}  < p^{\frac{1}{n+2}}.\eeq Consider
now the matrices $A_1=(1/r) B_1$ and $A_2=(1/r) B_2$. We claim
that $\rho(A_1^{\otimes k}+ A_2^{\otimes k})<1$ iff the instance
of 3-SAT is not satisfiable. Indeed, assume first that  the
instance of 3-SAT is not satisfiable. Then
$$\rho(B_1, B_2) \leq (p-1)^{1/(n+2)}.$$
Since $\rho(A_1^{\otimes k}+ A_2^{\otimes k})$ is an approximation
of relative accuracy $1/2^{2k}$ of $\rho(A_1, A_2)$, we have
$$ \rho(A_1^{\otimes k}+ A_2^{\otimes k}) \leq 2^{2k} \rho(A_1, A_2) \leq \frac{2^{2k}}{r} \rho(B_1, B_2) \leq \frac{2^{2k}}{r}
(p-1)^{1/(n+2)}$$ and by using the inequalities (\ref{e61}) we
conclude $\rho(A_1^{\otimes k}+ A_2^{\otimes k})<1$. In an
analogous way one see that, if the instance of 3-SAT is
satisfiable, then $\rho(A_1^{\otimes k}+ A_2^{\otimes k})>1$ and
the proof is therefore complete.

\qed

\section{Ellipsoid approximation}\label{sc-Ellips}
\label{s8} \setcounter{equation}{0}

Let us now compare our approximation of the joint spectral radius
with another approximation appearing in the literature. The
following approximation, called ellipsoid approximation, is
introduced in \cite{BNT}: \beq \label{prob-Hat} \hat \rho(A_1,
\dots , A_m) = \left[ \inf\limits_{X,\tau} \{ \tau: \; \tau X
\succeq A_i X A_i^T, \; i = 1, \dots, m,\; X \succ 0 \}
\right]^{1/2}. \eeq This expression corresponds to the best
ellipsoid norm for the set of matrices. In \cite{BNT} the
following inequalities are proved
$$
\ba{c}\frac{1}{\sqrt n}\ea \hat \rho(A_1, \dots , A_m) \leq
\rho(A_1, \dots , A_m) \leq \hat \rho(A_1, \dots , A_m).
$$
Thus, the ellipsoid approximation has a relative accuracy of
$1/\sqrt n$. This  accuracy decreases with the size of the
matrices but does not depend on the number of matrices in the set.
Using our results we can prove that the relative accuracy of the
ellipsoid approximation is in fact bounded by $1/ \sqrt m$ where
$m$ is the number of matrices.

\BT \label{t7} Let $\{ A_i: i=1, \ldots, m \}$ be a set of
matrices and define the ellipsoid approximation $\hat \rho$ by
(\ref{prob-Hat}). Then
$$
\ba{c}\frac{1}{\sqrt m}\ea \hat \rho(A_1, \dots , A_m) \leq
\rho(A_1, \dots , A_m) \leq \hat \rho(A_1, \dots , A_m).
$$
\ET

\proof It is easy to see that $\rho(A_1, \dots , A_m) \leq \hat
\rho(A_1, \dots , A_m)$. For proving the first inequality,  note
that the spectral radius of the linear operator
$$
B: X \rightarrow \sum\limits_{i=1}^m A_i X A_i^T: \quad {\cal S}
\to {\cal S}
$$
can be represented as follows: \beq\label{prob-B} \rho(B) =
\inf\limits_{X,\tau} \{ \tau:\; \tau X \succeq \sum\limits_{i=1}^m
A_i X A_i^T, \; X \succ 0 \}. \eeq If a pair $(X,\tau)$ is
feasible for the optimization problem in (\ref{prob-B}), then it
is also feasible for the optimization problem in (\ref{prob-Hat}).
Therefore
$$
\hat \rho(A_1, \dots , A_m) \leq \rho^{1/2}(B).
$$
Finally, by Theorem \ref{t5} we have
$$\ba{c}{1 \over \sqrt{m}}\ea
\rho^{1/2}(B) \leq \rho(A_1, \dots , A_m)$$ and so we get the
following bounds \beq\label{eq-AllBounds} \ba{c} {1 \over
\sqrt{m}} \hat \rho(A_1, \dots , A_m) \leq {1 \over \sqrt{m}}
\rho^{1/2}(B) \leq \rho(A_1, \dots , A_m) \leq \hat \rho(A_1,
\dots , A_m) \leq \rho^{1/2}(B). \ea \eeq Thus we see that the
ellipsoid estimate $\hat \rho(A_1, \dots , A_m)$ has relative
quality ${1 / \sqrt{m}}$. \qed

By using Theorem \ref{th-Kron} one can obtain a statement
analogous to that of Theorem \ref{t3}.

\begin{corollary}
\label{c51} Let $\{ A_i: i=1, \ldots, m \}$ be a set of matrices
and define the ellipsoid approximation $\hat \rho$ by
(\ref{prob-Hat}). Then \begin{equation} \label{e654} \ba{c}
\frac{1}{m^{1/(2k)}} \ea \hat \rho^{1/k}(A_1^{k}, \ldots, A_m^{k})
\leq \rho(A_1, \ldots, A_m)\leq \hat \rho^{1/k}(A_1^{k}, \ldots,
A_m^{k}).
\end{equation}
\end{corollary}

For $k-1$, this result is also proved in \cite{ando}. For
arbitrary $k$ a proof can be reconstructed from the result proved
in \cite{gur}. Notice however that approximation algorithms that
are directly based on (\ref{e654}) are not efficient because they
require to solve  linear matrix inequalities of very large
dimension. Approximations based on (\ref{e43}) are easier to
obtain  because there exists many efficient numerical algorithms
for approximating the spectral radius of a matrix.

\section{Discussion}\label{sc-Diss}
\label{s9} \setcounter{equation}{0}

The results presented in this paper introduce the possibility to
compute approximations for the joint spectral radius with
worst-case theoretical guarantees. Of course, the computational
complexity of these estimates needs further examination. Indeed,
the spectral structure of the matrix $A^{\otimes k}$ is very
simple and this simplicity must be inherited somehow by the matrix
$A_1^{\otimes k} + A_2^{\otimes k}$. It therefore appears to be an
interesting problem  to investigate the possibility of
constructing efficient procedures for finding the spectral radius
of this sum; of course there are theoretical limitations on what
can be achieved since we have shown that the problem of deciding
$\rho(A_1^{\otimes k}+A_2^{\otimes k}) \leq 1$ is NP-hard.\\

Another interesting issue is that of determining if worst-case
theoretical guarantees can be provided for the approximation
$$
\max\limits_{\sigma \in \{1,\dots,m\}^k} \rho (A_\sigma)^{1/k}
\leq \rho(A_1, \ldots, A_m).$$ In many numerical examples, this
approximation does in fact perform at least as well. \\

\end{document}